\documentclass{amsart}
\usepackage{amssymb}
\newcommand{\ZFCa}{{\operatorname{\mathsf {ZFC}}}}

\newcommand{\cf}{{\operatorname{\mathsf   {cf}}}}

\newcommand{\QED}{\hspace{0.1in} \square \vspace{0.1in}}
\newcommand{\forces}{\Vdash}
\newcommand{\F}{{\mathcal F}}
\newcommand{\N}{{\mathcal N}}
\newcommand{\M}{{\mathcal M}}
\newcommand{\V}{{\mathbf V}}
\newcommand{\HH}{{\mathbf H}}
\newcommand{\<}{\langle}
\renewcommand{\>}{\rangle}
\newcommand{\thinks}{\models}
\newcommand{\Proof}{{\sc Proof} \hspace{0.2in}}
\newcommand{\lft}[2]{\mathopen\ifcase#1{}\oo\or
                        \big#2\or\Big#2\else\oo\fi} 
\newcommand{\rgt}[2]{\mathclose\ifcase#1{}\oo\or
                        \big#2\or\Big#2\else\oo\fi} 
\newcommand{\PPP}{{\mathcal P}}
\theoremstyle{plain}
\newtheorem{theorem}{Theorem}[section]
\theoremstyle{plain}
\newtheorem{lemma}[theorem]{Lemma}

\begin{document}

\title{Remarks on the intersection of filters}
\author{Tomek Bartoszy\'{n}ski}
\address{Department of Mathematics\\
Boise State University\\
Boise, Idaho 83725, USA}
\email{{\tt tomek@@math.idbsu.edu}}

\thanks{Research partially supported by 
Idaho State Board of Education grant \#95--041 and NSF grant DMS 95-05375}
\maketitle

\begin{abstract}
We will show that that the existence of an uncountable family of
nonmeager filter whose intersection is meager is consistent with
${\mathbf{ MA}}(\text{Suslin})$.
  \end{abstract}
\section{Introduction}
Suppose that $\F$ is a filter on $\omega$. We identify elements of
$\F$ with their characteristic functions and think of $\F$ as a subset
of $2^\omega$. 
It is well known that if $\F$ is a nonprincipal filter (which we
assume to be always the case) then $\F$ is meager or $\F$ does not
have the property of Baire. Similarly, $\F$  has measure zero or is
nonmeasurable. Since an intersection of filters is again
a filter we want to know how many nonmeager filters one needs to
intersect to produce a meager filter. Define

\begin{itemize}
\item ${{\mathfrak  f}}_\M =\min\{|{\mathbf H}| : \forall \F \in
  {\mathbf H} \ \F$ is
a filter without the Baire property and  $\bigcap {\mathbf H}$ has
the Baire property $\},$
\item ${{\mathfrak  f}}_\N =\min\{|{\mathbf H}| : \forall \F \in  {\mathbf H} \ \F$ is
a nonmeasurable filter and 
$\bigcap {\mathbf H}$ has
measure zero$\}$.
\end{itemize}

\begin{theorem}[\cite{Tal80Com}]\label{2.1}
  $ {\mathfrak  t} \leq {\mathfrak  f}_\M$.
  In particular, 
  $\mathbf{MA}$ implies that the intersection of less than
  $2^{\boldsymbol\aleph_0}$ filters without the Baire property does
  not have the Baire property.
\end{theorem}

In \cite{Ple91Ide} Plewik found sharper estimates. In
particular, he  showed that $\mathfrak h \leq {\mathfrak f}_\M \leq
{\mathfrak d}$. Repicky improved the lower bound and showed that 
$\mathfrak g \leq {\mathfrak f}_\M$, where ${\mathfrak  g}=\min\{|{\mathbf H}| : \forall H \in  {\mathbf H} \
\text{ $H$ is 
groupwise dense and } \bigcap {\mathbf H} = \emptyset\}.$

Recall that a family $H \subseteq [\omega]^\omega $ is groupwise dense if
  \begin{enumerate}
  \item $\forall x \in  [\omega]^\omega \ [x]^\omega \cap H \neq
    \emptyset$,
  \item if $x \subseteq^\star y$ and $y \in H$, then $x \in H$, and 
  \item for every partition of $\omega$ into finite sets, $\{I_n : n
    \in \omega\}$, there is $a \in [\omega]^\omega$ such that
    $\bigcup_{n \in a} I_n \in H$. 
  \end{enumerate}

For measure the situation is quite different.
\begin{theorem}[Fremlin, \cite{BJbook}]\label{fremlinint}
  Assume $ {\mathbf {MA}}$.
Then there
  exists a family 
  $\left\{{\mathcal F}_{\xi} : \xi < 2^{\boldsymbol\aleph_{0}}\right\}$ of
  nonmeasurable filters such that $\bigcap_{\xi \in {\mathbf X}} {\mathcal F}_{\xi}$
  is a measurable filter for every uncountable set ${\mathbf X} \subseteq
  2^{\boldsymbol\aleph_{0}}$.  In particular, ${\mathbf {MA}}$ implies
  that there 
  exists a family of $\boldsymbol\aleph_{1}$ nonmeasurable filters
  with measurable 
  intersection.
\end{theorem}

The goal of this note is to show that the equality ${\mathfrak f}_\M=
\boldsymbol\aleph_1 $ is consistent with a relatively strong version
of Martin's Axiom.

Recall that a forcing notion $(\PPP, \leq)$ is Suslin if
\begin{enumerate}
  \item $\PPP$ is ccc,
  \item $\PPP$ is ${\boldsymbol \Sigma}^1_1$,
  \item relations $\leq, \perp$ are ${\boldsymbol \Sigma}^1_1$.
\end{enumerate}
Let ${\mathbf {MA}}(\text{Suslin})$ denote the Martin's Axiom for
Suslin partial orders.
It is well known that  ${\mathbf {MA}}(\text{Suslin})$ implies that
additivity of measure, $\mathfrak b$, etc. are all equal to
$2^{\boldsymbol\aleph_0}$.

\section{Consistency result}
The goal of this section is to show ${\mathbf {MA}}(\text{Suslin})$ is
consistent with
${\mathfrak  f}_\M=\boldsymbol\aleph_1$.

\begin{theorem}
There exists a model of $\V' \thinks \ZFCa$ and a family of filters
$\{\F_\alpha: \alpha < 2^{\boldsymbol\aleph_0}\} \in \V'$ such that:
\begin{enumerate}
  \item $\F_\alpha$ is not meager for each
    $\alpha$,
  \item $\bigcap_{\alpha \in \mathbf X} \F_\alpha$ is meager for every
    uncountable set $\mathbf X$,
  \item $\V' \thinks {\mathbf {MA}}(\text{Suslin}) +
    2^{\boldsymbol\aleph_0}=\boldsymbol\aleph_2$.
\end{enumerate}
\end{theorem}
\Proof
Let $\<{\mathcal P}_\alpha, \dot{{\mathcal Q}}_\alpha: \alpha<\omega_2\>$ be a finite support iteration such that 
\begin{enumerate}
  \item $\forces_\alpha \dot{{\mathcal Q}}_\alpha \text{ is Suslin,}$
  \item if $\alpha$ is a successor ordinal then $\dot{{\mathcal
        Q}}_\alpha$ adds a Cohen real.
\end{enumerate}
The second requirement is purely technical, its purpose it is to
simplify notation later on.

By careful bookkeeping we can ensure that $\V^{{\mathcal
    P}_{\omega_2}} \thinks {\mathbf {MA}}(\text{Suslin}) +
    2^{\boldsymbol\aleph_0}=\boldsymbol\aleph_2$.

Let ${\mathbf Z}=\{X_\alpha: \alpha <\omega_2\}$ be the  sequence of
Cohen reals added by $\dot{{\mathcal Q}}_\alpha$'s 
(represented as elements of $[\omega]^\omega$).

\begin{lemma}\label{lem1}
  ${\mathbf Z}$ is a generalized Luzin set. In
  particular, every subset of ${\mathbf Z}$ of size $
  \boldsymbol\aleph_2$ in nonmeager.
\end{lemma}
\Proof
Suppose that $A \subseteq [\omega]^\omega$ is a Borel meager set. Let
$\gamma$ be such that $A \in \V[{\mathcal P}_\gamma \cap G]$. Then $X_\alpha
\not\in A$ for $\alpha>\gamma$.~$\QED$

\begin{lemma}\label{cohens}
  If $\alpha_1< \alpha_2 < \dots < \alpha_n< \omega_2$ then
  $X_{\alpha_1} \cap X_{\alpha_2} \cap \dots \cap X_{\alpha_n}$ is a
  Cohen real over $\V[G \cap \PPP_{\alpha_1 -1}]$.
\end{lemma}
\Proof
Let 
$\phi:\lft1([\omega]^\omega\rgt1)^n \rightarrow [\omega]^\omega$ be
defined as $\varphi(X_1, \dots, X_n)=X_{1} \cap  \dots \cap X_{n}$.
Suppose that $A \subseteq [\omega]^\omega$ is a meager Borel set in
$\V[G \cap \PPP_{\alpha_1 -1}]$. 
Then $B_0 = \varphi^{-1}(A)$ is a meager set. 
Let $$C_0=\lft2\{X: \lft1\{\<X_2, \dots, X_n\>: \<X,X_2, \dots, X_n\>
\in B_0\rgt1\} 
\text{ is meager}\rgt2\}.$$ 
$C_0$ is a comeager set so $X_{\alpha_1} \in C_0$. 

Let 
$$B_1=\lft1\{\<X_2, \dots, X_n\>: \<X_{\alpha_1},X_2, \dots, X_n\>
\in B_0\rgt1\}.$$

$B_1$ is a meager set coded in $\V[G \cap \PPP_{\alpha_1}]
\subseteq \V[G \cap \PPP_{\alpha_2 -1}]$. 
We apply the construction above to get the set 
$$C_1=\lft1\{X: \{\<X_3, \dots, X_n\>: \<X,X_3, \dots, X_n\> \in B_1\}
\text{ is meager}\rgt1\},$$
and continue in this fashion.

It follows that 
$\<X_{\alpha_1}, X_{\alpha_2}, \dots, X_{\alpha_n}\> \not\in
\varphi^{-1}(A)$ which finishes the proof.~$\QED$

Let $\{Z_\alpha: \alpha<\omega_2\}$ be a partition of the set
$\{\alpha: \alpha=\beta+1, \beta<\omega_2\}$ into disjoint cofinal
sets.

Let $\F_\beta$ be a filter generated by sets $\{X_\alpha: \alpha
\in Z_\beta\}$. Since $X_\alpha$'s are Cohen reals each $\F_\beta$ is
indeed a filter.

Note that $\F_\beta \supseteq \{X_\alpha : \alpha \in Z_\beta\}$. Thus
it follows from \ref{lem1} that  $\F_\beta$ is not meager for every
$\beta<\omega_2$. 

\begin{theorem}
  $\bigcap_{\alpha\in \mathbf X} \F_\alpha$ is the filter of cofinite
  sets for every uncountable set $\mathbf X \subseteq \omega_2$.
\end{theorem}
\Proof
For simplicity assume that ${\mathbf X} = \omega_1$. The proof of the
general case is the same.

Suppose that $X \in \bigcap_{\xi<\omega_1} \F_\xi$ and $|\omega
\setminus X|=\boldsymbol\aleph_0$.

For each $\xi<\omega_1$ we can find $\alpha^\xi_1 < \alpha^\xi_2 <
\dots < \alpha^\xi_{n^\xi} \in Z_\xi$ such that 
$$X_{\alpha^\xi_1} \cap X_{\alpha^\xi_2} \cap 
\dots \cap X_{\alpha^\xi_{n^\xi}} \subseteq^\star X.$$

By passing to a subsequence we can assume that $n^\xi=n$ for all
$\xi<\omega_1$. 
Moreover, we can assume that $n$ is minimal, that is,
$$\forall \xi<\omega_1 \ X_{\alpha^\xi_1} \cap X_{\alpha^\xi_2} \cap 
\dots \cap X_{\alpha^\xi_{n-1}} \not \subseteq^\star X.$$

Let $\gamma$ be the least ordinal such that $\{\xi<\omega_1:
\alpha^\xi_n < \gamma\}$ is uncountable.
By the minimality of $n$, $\gamma$ is a limit ordinal and
$\cf(\gamma)=\boldsymbol\aleph_1$.
By passing to a subsequence again we can assume that
\begin{enumerate}
  \item $\alpha^\xi_n<\gamma$ for all $\xi<\omega_1$,
\item $\forall \delta < \gamma \ \exists \alpha \ \forall \xi>\alpha \
  \alpha^\xi_1 > \delta$ (since $Z_\alpha$'s are disjoint).
\end{enumerate}

We will show that $X \not\in \V[\PPP_\delta  \cap G]$ for all
$\delta<\omega_2$ and this contradiction will finish the proof.

For $\xi < \omega_1$ let 
$X_\xi= X_{\alpha^\xi_1} \cap X_{\alpha^\xi_2} \cap 
\dots \cap X_{\alpha^\xi_{n}}$. By the assumption $X_\xi
\subseteq^\star X$ for all $\xi< \omega_1$. 
\begin{lemma}
  $X \not\in \V[\PPP_\gamma  \cap G]$.
\end{lemma}
\Proof
Arguing as in \ref{lem1} and using  \ref{cohens} and condition (2)
above, we show that
$\{X_\xi: \xi < \omega_1\}$
is a Luzin set in 
$\V[\PPP_\gamma  \cap G]$. Since the set $\{Z \in [\omega]^\omega: Z
\subseteq^\star X\}$ is meager it follows that  $\{\xi: X_\xi
\subseteq^\star X\}$ is countable.~$\QED$

\begin{lemma}
  $X \not\in \V[\PPP_\delta   \cap G]$ for $\gamma <\delta<\omega_2$.
\end{lemma}
\Proof We will work in a model $\V'=\V[\PPP_\delta   \cap G]$. 
Suppose that the lemma is false. Let $\dot{X}$ be a
$\PPP_{\gamma,\omega_2}$-name for $X$. Let $M$ be a countable
elementary submodel of $\HH(\chi)$ containing $\dot{X}$ and
$\PPP_{\omega_2}$. 
Define a finite support iteration $\<\PPP_\alpha(M), \dot{{\mathcal
    Q}}_\alpha(M): \alpha<\omega_2\>$ as follows:

$$\forces_\alpha \dot{{\mathcal Q}}_\alpha = \left\{\begin{array}{ll}
\dot{{\mathcal Q}}_\alpha & \text{if $\alpha \in M$}\\
\emptyset & \text{if $\alpha \not\in M$}
\end{array}\right. \text{ for $\alpha<\omega_2$}. $$
Let ${\mathcal P}=\lim {\mathcal P}_\alpha(M)$.

${\mathcal P}$ is isomorphic to a countable iteration of Suslin
forcings. It may not be Suslin itself but it has enough absoluteness
properties to carry out the rest of the proof (see \cite{JuShSus88} or
\cite{BJbook} lemma 9.7.4). In particular, ${\mathcal P}$ has a definition that
can be coded as a real number. 

From Suslinness it follows  that ${\mathcal P} \lessdot {\mathcal
  P}_{\omega_2}$ and 
that $\dot{X}$ is a ${\mathcal P}$-name.

Let $N \prec \HH(\chi)$ be a countable model containing $M, \dot{X}$
and ${\mathcal P}$.
Since $\{X_\xi:\xi<\omega_1\}$ is a Luzin set in $\V'$ we can find
$\xi$ such that $Y=X_\xi$ is a Cohen real over $N$.
By the assumption $\forces_{{\mathcal P}} Y \subseteq^\star \dot{X}$.

By absoluteness, $N[Y][G \cap N[Y]] \thinks Y \subseteq^\star \dot{X}[G \cap N[Y]]$
and therefore  
$$N[Y] \thinks \text{``$\forces_{{\mathcal P}} Y \subseteq^\star
  \dot{X}$.''}$$
Represent Cohen algebra as $ \mathbf C = [\omega]^{< \omega}$ and let
$ \dot{Y}$ be the canonical name for a Cohen real.
There is a condition $p \in \mathbf C$ such that 
$$N \thinks p \forces_{\mathbf C} \text{``}\forces_{ {\mathcal P}} \dot{Y}
\subseteq^\star \dot{X}.\text{''}$$

Let $Y'= p \cup \lft1((\omega \setminus Y) \setminus
\max(p)\rgt1)$.
$Y'$ is also a a Cohen real over $N$ and since $p \subseteq Y'$ we get
that
$N[Y'] \thinks \text{``$\forces_{{\mathcal P}} Y' \subseteq^\star
  \dot{X}$.''}$
It follows that 
$$N[Y'][G \cap N[Y']] \thinks Y' \subseteq^\star \dot{X}[G
\cap N[Y']].$$ 
Note that $ \dot{X}[G] = \dot{X}[G
\cap N[Y']]=\dot{X}[G \cap N[Y]]$.
Thus $\V[G] \thinks Y \cup Y' \subseteq^\star \dot{X}[G]$ which means that $
\dot{X}[G]$ is cofinite. Contradiction.~$\QED$

\ifx\undefined\bysame
\newcommand{\bysame}{\leavevmode\hbox to3em{\hrulefill}\,}
\fi

\end{document}